\documentclass[letter, 11pt]{article}
\usepackage{amssymb, amsthm, graphicx}
\usepackage{mathrsfs}
\usepackage{booktabs,threeparttable}
\usepackage[top=1in, bottom=1in, left=1in, right=1in]{geometry}
\usepackage{graphicx}

\newtheorem{proposition}{Proposition}[section]
\newtheorem{theorem}{Theorem}[section]
\newtheorem{definition}{Definition}[section]

\newtheorem{lemma}{Lemma}[section]

\begin{document}

\title{A New Condition for the Existence of Optimal Stationary Policies in Denumerable State Average Cost Continuous Time Markov Decision Processes with
Unbounded Cost and Transition Rates}

\author{Ping Cao, Jingui Xie, \\ School of Management, University of Science and Technology of China, \\ pcao@ustc.edu.cn, xiej@ustc.edu.cn}

\maketitle

{\bf Abstract} This paper presents a new condition for the existence of optimal stationary policies in average-cost continuous-time Markov decision processes with unbounded cost and transition rates, arising from controlled queueing systems. This condition is closely related to the stability of queueing systems. It suggests that the proof of the stability can be exploited to verify the existence of an optimal stationary policy. This new condition is easier to verify than existing conditions. Moreover, several conditions are provided which suffice for the average-cost optimality equality to hold.
\\
\textit{Keywords}: Markov decision processes; Average-cost criterion; Unbounded transition rates; Optimal stationary policy







\section{Introduction}

Queueing systems have wide applications in computer communication networks, manufacturing processes and customer service platforms~\cite{cmp}. There exists a lot of literature studying issues such as system stability, cost and performance analysis under a given service principle, e.g.,~\cite{awz,xhz08}. In order to cut down the operational cost and better serve customers, the queueing models should be controlled in a way such that the operational cost is minimized, e.g.,~\cite{h84}. A lot of controlled queueing models can be analyzed as continuous-time Markov decision processes (CTMDP)~\cite{ws87}. The buffer of the queue model is often unlimited, and transition rates might be dependent on the system state. Therefore, the corresponding CTMDP often has denumerable states and the transition rates are unbounded. Moreover, the state-dependent cost rates are also unbounded.  A question naturally arises that whether an optimal stationary policy for such a CTMDP exists or not. For the discounted cost CTMDP, it is often relatively easy to verify whether an optimal stationary policy exists~\cite{guoher,guoz02a}. However, for the average-cost CTMDP, more conditions should be imposed to ensure the existence of an optimal stationary policy~\cite{guoher,guoz02b}. This paper provides a new condition under which an average-cost optimal stationary policy exists, which is different from existing conditions.

In 2002, \cite{guoz02b} presents a set of conditions under which an average-cost optimal stationary policy exists. Their conditions requires constructing a series of functions which satisfy their proposed assumptions. It is not straightforward to construct these functions for most problems we encounter. Later in 2009, \cite{guoher} gives a sufficient condition for the existence of optimal stationary policies, which also requires finding a function satisfying several conditions. However, this function is often problem specific. Without adequate research of the specific CTMDP, it is not easy to find an appropriate function. It will be valuable if we can find a way to bypass seeking for such a function.
Focused on discrete-time Markov decision processes (DTMDP), \cite{sen99} gives several conditions under which an average cost optimal stationary policy exists.
\cite{sen99} has also mentioned that CTMDP can be transformed into DTMDP if the transition rate is uniformly bounded by employing uniformization method. However, the uniformization method cannot be applied if the transition rate is unbounded. Therefore, it is quite necessary to analyze the CTMDP with unbounded transition rates separately.

In this paper, a new condition is provided to ensure the existence of an average-cost optimal stationary policy for denumerable state CTMDP with unbounded cost and transition rates. This condition concerns that whether the expected time and expected cost of a first passage from any state to a given state is finite or not under a given controlled policy. The former is related to the stability of the queueing system, while the latter may be seen as a generalized stability if we notice that the expected cost is equal to the expected time if the cost rate is 1. A lot of literature has focused on discussing the stability of the queueing system, e.g.,~\cite{t75,xhz08}. Their results can help us verify whether an average-cost optimal stationary policy exists or not.

This paper is organized as follows. In Section 2, we introduce CTMDP and present our main result. Section 3 gives the proof of the main result.
In Section 4, we give conditions under which the average-cost optimality inequality (ACOI) becomes equality.

\section{Model and Main Result}

Consider a continuous-time Markov decision process $\{x(t): t\geq 0\}$ consisting of four-element tuple $\{S, (A(i), i\in S), q(j|i,a), c(i,a)\}$:
\begin{enumerate}
  \item The state space $S$ is denumerable;
  \item Each action space $A(i)$ is a subset of the finite action space $A$;
  \item The transition rate $q(j|i,a)$ satisfies $q(j|i,a)\geq 0$, $\forall\ i\neq j$, $i, j\in S$, $a\in A(i)$ and $\sum_{j\in S}q(j|i,a)=0$, $\forall\ i\in S, a\in A(i)$.
  \item The cost rate function $c(i,a)\geq 0$, $\forall\ i\in S, a\in A(i)$.
\end{enumerate}

Let $\Pi$ be the set of all randomized Markov policies and $F$ the set of all stationary policies~\cite{guoher}. Given $\pi=(\pi_t)\in\Pi$ and the discount factor $\alpha>0$, we define the expected discounted cost function (with initial state $ i $)
\begin{equation}
J_{\alpha}(i,\pi)=\int_{0}^{\infty}e^{-\alpha t}E_{i}^{\pi}[c(x(t),\pi_t)]dt, \forall\ i\in S, \pi\in\Pi, \label{eq_dis}
\end{equation}
and the corresponding optimal discounted cost function
$J_{\alpha}^{*}(i)=\inf_{\pi\in\Pi}J_{\alpha}(i,\pi), \forall\ i\in S,$
where $c(i,\pi_t)$ is the expected cost rate at state $i$ using policy $\pi_t$ at time $t$, which is defined as
$c(i,\pi_t)=\int_{A(i)}c(i,a)\pi_t(da|i)$.

It is straightforward to show that for any stationary policy $f\in F$, we have
\begin{equation}
\alpha J_{\alpha}(i,f)=c(i,f(i))+\sum_{j\in S}J_{\alpha}(j)q(j|i,f(i)). \label{eq_dis_dp_f}
\end{equation}

Since $A(i)$ is finite, \cite{guoher} states that $J_{\alpha}^{*}(i)$ is well defined and satisfies the discounted-cost optimality equation
\begin{equation}
\alpha J_{\alpha}^{*}(i)=\min_{a\in A(i)}\left\{c(i,a)+\sum_{j\in S}J_{\alpha}^{*}(j)q(j|i,a)\right\}. \label{eq_dis_dp}
\end{equation}

Moreover, we define the long run expected average cost function
\begin{equation}
J_c(i,\pi)=\lim\sup_{T\rightarrow\infty}\frac{1}{T}\int_0^TE_{i}^{\pi}[c(x(t),\pi_t)]dt, \forall\ i\in S, \pi\in\Pi, \label{eq_ave}
\end{equation}
and the corresponding optimal average cost function
$J_{c}^{*}(i)=\inf_{\pi\in\Pi}J_{c}(i,\pi), \forall\ i\in S.$

One of the most important questions is whether an average-cost optimal stationary policy exists for the CTMDP. Before we state our main results, we propose the following definition~\cite{sen99}.

{\definition Let $d$ be a (randomized) stationary policy. Then $d$ is a $i_0$-standard policy if the Markov process induced by $d$, $\{x^d(t):t\geq 0\}$ satisfies that for any $i\in S$, the expected time $m_{i,i_0}(d)$ of a first passage from $i$ to $i_0$ (during which  at least one transition occurs) is finite and the expected cost $c_{i,i_0}(d)$ of a first passage from $i$ to $i_0$ (during which at least one transition occurs) is finite.}

{\bf Remark 1:} Note that $x(t)=x(t+)$, a.e.. Thus, if we define the first passage time $\tau_{i,i_0}$ as $\tau_{i,i_0}=\inf\{t>0: x(t)=i_0|x(0)=i\}$, then $\tau(i,i_0)=0$ a.e. if $i=i_0$. Hence we impose additional constraint that at least one transition occurs on the definition of the first passage time.

{\bf Remark 2:} If the cost rate function is bounded, then $m_{i,i_0}(d)<\infty$ can implies $c_{i,i_0}(d)<\infty$. In this case, $d$ is a $i_0$-standard policy if the Markov process induced by $d$ is ergodic (i.e., irreducible and positive recurrent).

The following lemma is extensively used for analysing the stability of a queueing systems, of which the proof is omitted for brevity.

{\lemma\label{lem_1} Assume that $m_{i,i_0}<\infty$, $\forall\ i\in S$. Assume that there exists a (finite) nonnegative function $r$ on $S$ and a finite subset $H^*$ containing $i_0$ such that
\begin{equation}
\sum_{j}q(j|i)r(j)<\infty, i\in H^*, \label{eq_lem_1}
\end{equation}
and
\begin{equation}
c(i)+\sum_{j}q(j|i)r(j)\leq 0, i\notin H^*. \label{eq_lem_2}
\end{equation}
Then there exists a (finite) nonnegative constant $F$ such that $c_{i,i_0}\leq r(i)-r(i_0)+Fm_{i,i_0}$, $\forall\ i\neq i_0$. Especially, if $H^*=\{i_0\}$, then $c_{i,i_0}\leq r(i)$, $\forall\ i\neq i_0$.
}

Let $ S = \{0,1,2,\dots\} $. Now we propose our main result.

{\theorem\label{the_2} Assume that $J_{\alpha}^*(i)$ is increasing in $i$ for $\alpha>0$. If there exists a $0$-standard policy $d$, then there exists a constant $g^*\geq 0$, a stationary policy $f^*$, and a real-valued function $h^*$ (which is increasing in $i$) such that:

(i) There exists a sequence $\{\alpha_n, n\geq 1\}$ tending to zero (as $n\rightarrow\infty$) such that $\forall\ i\in S$,
\begin{equation}
f^*(i)=\lim_{k\rightarrow\infty}f_{\alpha_k}^*(i), g^*=\lim_{k\rightarrow\infty}\alpha_k J_{\alpha_k}^*(0),\label{eq_limit_fg}
\end{equation}
and
\begin{equation}
h^*(i)=\lim_{k\rightarrow\infty}h_{\alpha_k}(i), \label{eq_limit_h}
\end{equation}
where $h_{\alpha}(i):=J_{\alpha}^*(i)-J_{\alpha}^*(0)$.

(ii) $(g^*, f^*, h^*)$ satisfy the following average-cost optimality inequality (ACOI):
\begin{eqnarray}
g^* &\geq& c(i,f^*)+\sum_{j\in S}h^*(j)q(j|i,f^*)  \label{eq_acoe}
\\ &= & \min_{a\in A(i)}\left\{c(i,a)+\sum_{j\in S}h^*(j)q(j|i,a)\right\}, \forall i\in S, \nonumber
\end{eqnarray}
and $f^*$ is an average-cost optimal stationary policy.
}

{\bf Remark 1:} The above result still holds when the state is a vector rather than a scalar.

{\bf Remark 2:} The monotonicity of the discounted value function $J_{\alpha}^*(i)$ is often satisfied, e.g.,
in queueing systems more customers staying in the queue implies more waiting.

{\bf Remark 3:} The $0$-standard policy $d$ is not required to be optimal. It can be any policy which is easy to be constructed and analyzed.

{\bf Remark 4:} This theorem closely relates the existence of an average-cost optimal stationary policy to the stability of the queueing system under a given service policy. The queueing system is called to be stable under a given service policy if the induced Markov process is ergodic (irreducible and positive recurrent). Positive recurrence implies that for any $i\in S$, the expected time $m_{i,0}$ of a first passage from $i$ to $0$ is finite. To prove (positive recurrence) the finiteness of the expected time $m_{i,0}$, a Lyapunov function $r(\cdot)$ might be constructed in order to apply Lemma~\ref{lem_1} with $c(i)=1$ and $H^*=\{0\}$ (noticing that the expected cost is expected time if the cost rate is 1). This method can also be found in Theorem 1.18 in~\cite{c91}. In many situations, with slightly modification of the Lyapunov function $r(\cdot)$ constructed for proving $m_{i,0}<\infty$, another Lyapunov function can be constructed to satisfy (\ref{eq_lem_2}) and thus the finiteness of the expected cost $c_{i,0}$ can be proved. That is to say, the discussion of stability of the queueing system can help prove the existence of  an average-cost optimal stationary policy.

\section{Proof of Theorem~\ref{the_2}}

\cite{guoher} proposes the following assumptions to ensure the existence of an average-cost optimal stationary policy, which can be seen as a continuous-time counterpart of (SEN) assumptions proposed in~\cite{sen99}.

{\bf Assumptions A:} For some decreasing sequence $\{\alpha_n, n\geq 1\}$ tending to zero (as $n\rightarrow\infty$) and some state $i_0\in S$,

(A1) $\alpha_n J_{\alpha_n}^*(i_0)$ is bounded in $n$.

(A2) There exists a nonnegative (finite) function $H$ such that $h_{\alpha_n}(i)\leq H(i)$, $\forall\ i\in S, n\geq 1$, where $h_{\alpha}(i)=J_{\alpha}^*(i)-J_{\alpha}^*(i_0)$.

(A3) There exists a nonnegative constant $L$ such that $-L\leq  h_{\alpha_n}(i)$, $\forall\ i\in S, n\geq 1$.



Before proving Theorem~\ref{the_2}, we give some results of the Markov process $\{x(t): t\geq 0\}$.
Let
$J_{i}(t)=\frac{1}{t}E\left[\int_0^t c(x(s))ds|x(0)=i\right], \forall\ i\in S.$
We have the following result.

{\proposition\label{pro_5} Let $R$ be a positive recurrent class.

(i) For $i\in R$, $\lim_{t\rightarrow\infty}J_i(t)$ exists and equals the (finite or infinite) constant $J_R=: \sum_{j\in R}\pi_j c(j)$, where $\pi_j$ is the steady sate probability of being in state $j$.

(ii) For $i\in R$, we have $J_R=c_{i,i}/m_{i,i}$.

(iii) $J_R=\sum_{j\in R}\pi_j E[c(x(t))|x(0)=j]$, $\forall\ t\geq 0$.

}

\begin{proof}
Let $e_{i,j}$ be the expected time of visits to $j$ during a first passage from $i$ to $i$. Then $\pi_j=e_{i,j}/m_{i,i}$. Therefore, $J_R=\sum_{j\in R}\pi_j c(j)=\sum_{j\in R}c(j)e_{i,j}/m_{i,i}=c_{i,i}/m_{i,i}$ and thus (ii) holds.

Note that $J_{i}(t)=\sum_{j}c(j)E[\int_0^ t 1(x(s)=j)ds|x(0)=i]t^{-1}$ and $\lim_{t\rightarrow\infty}E[\int_0^ t 1(x(s)=j)ds|x(0)=i]t^{-1}=\pi_j$. By Fatou lemma, it follows that $\lim\inf_{t\rightarrow\infty}J_i(t)\geq J_R$. Thus, if $J_R=\infty$, the limit exists and equals $\infty$, $\forall\ i\in R$.
If $J_R<\infty$, then (i) follows from the renewal reward theorem (See~\cite{ross96}).

Next we prove (iii). Note that $E[c(x(t))|x(0)=j]=\sum_{k\in S}p(j,k,t)c(k)$, where $p(j,k,t)=P[x(t)=k|x(0)=j]$. Since $\sum_{j\in R}\pi_j p(j,k,t)=\sum_{j\in S}\pi_j p(j,k,t)=\pi_k$ (noting that $\pi_i=0$ for $i\in S-R$), we have
\begin{eqnarray*}
&&\sum_{j\in R}\pi_j E[c(x(t))|x(0)=j]=\sum_{j\in R}\pi_j\sum_{k\in S}p(j,k,t)c(k)
\\&=& \sum_{k\in S}c(k)\sum_{j\in R}\pi_j p(j,k,t)= \sum_{k\in S}c(k)\pi_k=\sum_{k\in R}c(k)\pi_k=J_R,
\end{eqnarray*}
where the interchange of the order of summation is valid as all terms are nonnegative.
\end{proof}

{\proposition Suppose that $d$ is a $i_0$-standard policy with positive recurrent class $R$. Let $J_R(d)$ and $\pi_i(d)$ be defined as in Proposition~\ref{pro_5}, then
\begin{equation}
J_R(d)=\alpha\sum_{i\in R}\pi_i(d)J_{\alpha}(i,d), \forall\ \alpha>0. \label{eq_standard}
\end{equation}
}

\begin{proof}
It follows from (\ref{eq_dis}) and Proposition~\ref{pro_5}(iii) that
\begin{eqnarray*}
&&\alpha\sum_{i\in R}\pi_i(d)J_{\alpha}(i,d)=\alpha\sum_{i\in R}\pi_i(d)\int_{0}^{\infty}e^{-\alpha t}E^{d}[c(x(t),d)|x(0)=i]dt
\\&=&\alpha\int_{0}^{\infty}e^{-\alpha t}\left[\sum_{i\in R}\pi_i(d)E^{d}[c(x(t),d)|x(0)=i]\right]dt=J_R(d),
\end{eqnarray*}
where the interchange of the summation and integration is valid as all terms are nonnegative.
\end{proof}

{\proposition\label{pro_3} Assume that $J_{\alpha}^*(i_0)<\infty$, for some $\alpha>0$. Given $i\neq i_0$, assume that there exists a policy $\theta_i$ such that both the expected time and expected cost of a first passage from $i$ to $i_0$ are finite. Then $h_{\alpha}(i)\leq c_{i,i_0}(\theta_i)$, and hence (A2) holds for $i_0$ with $H(i)=c_{i,i_0}(\theta_i)$.}

\begin{proof}
If the process begins in state $i\neq i_0$ and follows policy $\theta_i$, it will reach state $i_0$ at some time in the future, which is denoted by $T$. Let the policy $\psi$ follow $\theta_i$ until $i_0$ is reached, then follow an $\alpha$ discounted optimal policy $f_{\alpha}$.

Then we have
\begin{eqnarray}
J_{\alpha}^*(i)&\leq & J_{\alpha}(i,\psi) \nonumber
\\&=&E^{\psi}\left[\int_0^Te^{-\alpha t}c(x(t), a(t))dt|x(0)=i\right]+E^{\psi}\left[e^{-\alpha T}|x(0)=i\right]J_{\alpha}^*(i_0) \nonumber
\\&\leq & E^{\psi}\left[\int_0^Tc(x(t), a(t))dt|x(0)=i\right]+J_{\alpha}^*(i_0) \nonumber
\\&\leq & c_{i,i^0}(\theta_i)+J_{\alpha}^*(i_0). \label{eq_com}
\end{eqnarray}
The result follows by subtracting $J_{\alpha}^*(i_0)$ from both sides.
\end{proof}

{\bf Remark:} Proposition~\ref{pro_3} gives a way to construct a function $H(i)$. From the remark below Proposition~\ref{pro_2}, it is known that
$c_{i,i_0}(\theta_i)$ is a quite good choice for $H(i)$.

{\bf Proof of Theorem~\ref{the_2}: } We only need to prove that (A1-3) hold under conditions in Theorem~\ref{the_2}. Let $i_0=0$. It follow from (\ref{eq_standard}) that $J_R(d)\geq \alpha\pi_{0}(d)J_{\alpha}(0,d)\geq \alpha\pi_{0}(d)J_{\alpha}^*(0)$. Hence $\alpha J_{\alpha}^*(0)\leq J_R(d)/\pi_{0}(d)=c_{0,0}(d)$. Therefore, (A1) holds. From Proposition~\ref{pro_3} we know that (A2) holds with $H(i)=c_{i,0}(d)$ for $i\neq 0$ and $H(0)=0$. Since $J_{\alpha}^*(i)$ is increasing in $i$, it follows that $h_{\alpha}(i)\geq 0$, and hence (A3) holds with $L=0$. It follows from (\ref{eq_limit_h}) and the fact that $h_{\alpha}(i)$ is increasing in $i$ that $h^*(i)$ is increasing in $i$. $\hfill\square$

\section{Sufficient Conditions for ACOE to Hold}

Proposition 5.11, \cite{guoher}  has given an example to demonstrate that (SEN-C) is not sufficient to claim that the average-cost optimality equality (ACOE) holds, i.e., ACOI might be strict. \cite{guoher} gives one condition under which  the ACOE holds, i.e., the inequality in (\ref{eq_acoe}) is in fact equality. However, in many situations it is hard to verify this condition and even in some cases it fails to hold due to improper choice of the function $H(i)$.

In this section, we give conditions under which the ACOE holds.
We first develop some notations. Let $\mathfrak{R}(i,G)$ be the class of policies $\theta$ satisfying
$$P_{\theta}(x(t)\in G \mbox{ for some } t>0, \mbox{at least one transition occurs}|x(0)=i)=1,$$
 and the expected time $m_{i,G}(\theta)$ of a first passage from $i$ to $G$ (during which  at least one transition occurs) is finite. Let $\mathfrak{R}^*(i,G)$ be the class of policies $\theta\in\mathfrak{R}(i,G)$ such that the expected cost $c_{i,G}(\theta)$ of a first passage from $i$ to $G$ (during which  at least one transition occurs) is finite. If $G=\{x\}$, then $\mathfrak{R}(i,G)$ is denoted by $\mathfrak{R}(i,x)$ (respectively, $\mathfrak{R}^*(i,G)$ by $\mathfrak{R}^*(i,x)$).

%
%

{\proposition\label{pro_2} Assume that the Assumptions (A1-3) hold, and for some state $i$ and nonempty set $G$, there exists a policy $\theta\in\mathfrak{R}(i,G)$ such that $\sum_{j\in G}H(j)P_{\theta}(x(T)=j)<\infty$, where $T$ is the first passage time from $i$ to $G$ and $H$ is the function from (A2). Then for any limit function $h^*$, we have
\begin{equation}
h^*(i)\leq c_{i,G}(\theta)-g^*m_{i,G}(\theta)+E_{\theta}[h^*(x(T))|x(0)=i]. \label{eq_recur}
\end{equation}
}

\begin{proof}
In a derivation very similar to that in (\ref{eq_com}), we have
$$J_{\alpha}^*(i)\leq c_{i,G}(\theta)+E_{\theta}[e^{-\alpha T}J_{\alpha}^*(x(T))|x(0)=i],$$
which can be written as
\begin{equation}
h_{\alpha}(i)\leq c_{i,G}(\theta)-\alpha J_{\alpha}^*(i_0)\left(\frac{1-E_{\theta}[e^{-\alpha T}|x(0)=i]}{\alpha}\right)+E_{\theta}[e^{-\alpha T}h_{\alpha}(x(T))|x(0)=i]. \label{eq_com_2}
\end{equation}
Note that
\begin{equation}
\frac{1-E_{\theta}[e^{-\alpha T}|x(0)=i]}{\alpha}=E_{\theta}\left[\int_0^Te^{-\alpha s}ds|x(0)=i\right]. \label{eq_com_1}
\end{equation}

The term $\int_0^Te^{-\alpha s}ds$ is decreasing in $\alpha$. It follows from monotone convergence theorem that the limit of the left side of (\ref{eq_com_1}) exists and equals to $E_{\theta}[T|x(0)=i]=m_{i,G}(\theta)$.

Choose a discount factor sequence $\{\alpha_n, n\geq 1\}$ tending to zero such that (\ref{eq_limit_fg}) and (\ref{eq_limit_h}) hold. Taking the limit of both sides of (\ref{eq_com_2}) as $\alpha_n\rightarrow 0^+$ yields
$$
h^*(i)\leq c_{i,G}(\theta)-g^*m_{i,G}(\theta)+\lim_{n\rightarrow\infty}E_{\theta}[e^{-\alpha_n T}h_{\alpha_n}(x(T))|x(0)=i].
$$

Note that $e^{-\alpha_n T}h_{\alpha_n}(x(T))$ converges to $h^*(x(T))$ as $n\rightarrow\infty$. Since $e^{-\alpha_n T}h_{\alpha_n}(x(T))$ is bounded by $\max(L, H(x(T)))$ from (A2) and (A3), and $E_{\theta}[\max(L, H(x(T)))]\leq L+ E_{\theta}H(x(T))=L+\sum_{j\in G}H(j)P_{\theta}(x(T)=j)<\infty$, by dominated convergence theorem it is known that
$$\lim_{n\rightarrow\infty}E_{\theta}[e^{-\alpha_n T}h_{\alpha_n}(x(T))|x(0)=i]=E_{\theta}[h^*(x(T))|x(0)=i].$$

Therefore, (\ref{eq_recur}) holds.
\end{proof}

Now we give sufficient conditions under which the ACOE holds.

{\theorem\label{the_1} Assume that the Assumptions (A1-3) hold, and let $e$ be a stationary policy realizing the minimum in the ACOI. Define the nonnegative discrepancy function $\Phi$ to satisfy
\begin{equation}
g^*=c(i,e)+\Phi(i)+\sum_{j\in S}q(j|i,e)h^*(j), i\in S. \label{eq_ave_op_dis}
\end{equation}

Then $\Phi(i)=0$, and hence the ACOE holds at the particular state $i$ under any of the following conditions:

(i) There exists a nonempty set $G$ such that $e$ satisfies $e\in \mathfrak{R}(i,G)$ and $\sum_{j\in G}H(j)P_{\theta}(x(T)=j)<\infty$, where $T$ is the first passage time from $i$ to $G$.

(ii) $e\in\mathfrak{R}(i,i_0)$.

(iii) The Markov process induced by $e$ is positive recurrent at $i$.

(iv) $\sum_{j\in S}|q(j|i,a)|H(j)<\infty$ for $a\in A(i)$. This conditions typically hold when the jump size at each state $i$ is bounded and thus there are finite number of $j$ such that $q(j|i,a)>0$ for each $i\in S$.
}

\begin{proof}

To prove equality under (i), let the process operate under $e$, and suppress the initial state $i$.
Since the first passage time from $i$ to $G$, $T$, is a stopping time such that $E_e[T]=m_{i,G}(e)<\infty$ as $e\in \mathfrak{R}(i,G)$, it follows from Dynkin's formula (see~\cite{os05}) that
$$E^e[h^*(x(T))]=h^*(i)+E^e\left[\int_0^T \sum_{j\in S}q(j|x(s),e)h^*(j)ds\right].$$

From (\ref{eq_ave_op_dis}) it is known that
$$
E_e[h^*(x(T))]=h^*(i)+E^e\left[\int_0^T (g^*-c(x(s),e)-\Phi(x(s))ds\right],
$$
and thus
\begin{equation}
c_{i,G}(e)-g^*m_{i,G}(e)+E^e\left[\int_0^T\Phi(x(s))ds\right]+E^e[h^*(x(T))]=h^*(i), \label{eq_ave_op_dis_2}
\end{equation}
which implies that $c_{i,G}(e)<\infty$, and hence $e\in \mathfrak{R}^*(i,G)$. Therefore, we can apply Proposition~\ref{pro_2}, which yields
$$c_{i,G}(e)-g^*m_{i,G}(e)+E^e[h^*(x(T))]\geq h^*(i).$$
Comparing the above equation with (\ref{eq_ave_op_dis_2}) and keeping in mind that $\Phi$ is nonnegative, we know that $\Phi=0$ during the first passage from $i$ to $G$. Specially, we have $\Phi(i)=0$ and thus the ACOE holds at state $i$.

(ii) follows from (i) by choosing $G=\{i_0\}$ and the fact $h^*(i_0)=0$.

(iii) follows from (i) by noting that if the Markov process induced by $e$ is positive recurrent at $i$, then $e\in \mathfrak{R}(i,i)$.

(iv) follows from the same argument in Theorem 5.9 in~\cite{guoher}.
\end{proof}

{\bf Remark:} If starting from an arbitrary initial state $i$, in a finite expected amount of time the Markov process induced by $e$ reaches a finite set $G$, then the ACOE holds.

\section{A Queueing Example}

{\bf Example 1.} A single-server, 2-buffer queueing model. Consider a server serving two types of customers: type 1 and type 2 customers. Type 1 and 2 customers form queue 1 and queue 2, respectively. Type 1 and 2 customers arrive according to two independent Poisson processes with parameter $\lambda_1$ and $\lambda_2$, respectively. Buffers of both queues are assumed to be infinitely large. The service times of type 1 and 2 customers are exponentially distributed with parameters $\mu_1$ and $\mu_2$, respectively. While waiting in queue, a type 1 customer may change to a type 2 customer after a random time T, which is exponentially distributed with parameter $\lambda_T$. The holding cost of a customer in queue 1 and 2 per unit time is $h_1$ and $h_2$, respectively. When a type 1 customer upgrades, the cost of transferring from queue 1 to queue 2 is $c$ per unit. The server should decide which buffer to serve to minimize the average cost.

The state can be denoted by $\mathbf{q}=(q_1,q_2)$, where $q_i$ is the length of queue $i$, $i=1,2$. For each state $\mathbf{q}$, we have the corresponding action set
$$A(\mathbf{q})=\left\{\begin{array}{ll}\{0\}, & \mbox{ if } \mathbf{q}=(0,0),
\\ \{1\}, & \mbox{ if } \mathbf{q}=(q_1,0), q_1>0,
\\ \{2\}, & \mbox{ if } \mathbf{q}=(0,q_2), q_2>0,
\\ \{1,2\}, & \mbox{ otherwise. }
\end{array}
\right.
$$

And the corresponding transition rate is
$$q(\mathbf{q}'|\mathbf{q}, 1)=\left\{\begin{array}{ll}\mu_1, & \mbox{ if } \mathbf{q}'=(q_1-1,q_2),
\\ \lambda_1, & \mbox{ if } \mathbf{q}'=(q_1+1,q_2),
\\ \lambda_2, & \mbox{ if } \mathbf{q}'=(q_1,q_2+1),
\\ q_1\lambda_T, & \mbox{ if } \mathbf{q}'=(q_1-1,q_2+1),
\\ -(\mu_1+\lambda_1+\lambda_2+q_1\lambda_T), & \mbox{ if } \mathbf{q}'=(q_1,q_2),
\\ 0, & \mbox{ otherwise;}
\end{array} \mbox{ for } q_1\geq 1,
\right.
$$
and
$$q(\mathbf{q}'|\mathbf{q}, 2)=\left\{\begin{array}{ll}\mu_2, & \mbox{ if } \mathbf{q}'=(q_1,q_2-1),
\\ \lambda_1, & \mbox{ if } \mathbf{q}'=(q_1+1,q_2),
\\ \lambda_2, & \mbox{ if } \mathbf{q}'=(q_1,q_2+1),
\\ q_1\lambda_T, & \mbox{ if } \mathbf{q}'=(q_1-1,q_2+1),
\\ -(\mu_2+\lambda_1+\lambda_2+q_1\lambda_T), & \mbox{ if } \mathbf{q}'=(q_1,q_2),
\\ 0, & \mbox{ otherwise;}
\end{array} \mbox{ for } q_1\geq 1, q_2\geq 1,
\right.
$$

$$q(\mathbf{q}'|\mathbf{q}, 2)=\left\{\begin{array}{ll}\mu_2, & \mbox{ if } \mathbf{q}'=(q_1,q_2-1),
\\ \lambda_1, & \mbox{ if } \mathbf{q}'=(q_1+1,q_2),
\\ \lambda_2, & \mbox{ if } \mathbf{q}'=(q_1,q_2+1),
\\ -(\mu_2+\lambda_1+\lambda_2), & \mbox{ if } \mathbf{q}'=(q_1,q_2),
\\ 0, & \mbox{ otherwise;}
\end{array} \mbox{ for } q_1=0, q_2\geq 1,
\right.
$$

Moreover, $q(\mathbf{q}'|\mathbf{q}, 0)=0$.

Let $x(t)=(x_1(t),x_2(t))$ be state at time $t$, and $T_{r}^\pi(t)$ be the total number of transferred customers till time $t$ under policy $\pi$, where $x_i(t)$ is the length of queue $i$ at time $t$, $i=1,2$. The expected discounted cost function under policy $\pi$ in this example can be formulated as
\begin{eqnarray*}
J_{\alpha}(\mathbf{q},\pi)&=&E_{\mathbf{q}}^{\pi}\left[\int_{0}^{\infty}e^{-\alpha t}(dT_{r}^\pi(t)+(h_1x_1(t)+h_2x_2(t))dt)\right]
\\&=& \int_{0}^{\infty}e^{-\alpha t}E_{\mathbf{q}}^{\pi}(h_1x_1(t)+h_2x_2(t)+c\lambda_T x_1(t))dt.
\end{eqnarray*}
Hence, the cost rate function is $c(\mathbf{q},1)=c(\mathbf{q},2)=c(\mathbf{q})=h_1q_1+h_2q_2+cq_1\lambda_T$.

%
%
%
%

We have the following result for Example 1.

{\proposition Suppose that $\lambda_1+ \lambda_2< \min(\mu_1,\mu_2)$. There exists an average-cost optimal stationary policy for Example 1 and the ACOE holds.
}

\begin{proof}
We apply Theorem~\ref{the_2} by proving that

(i) $J_{\alpha}^{*}(\mathbf{q})$ is increasing in $\mathbf{q}$;

(ii) The priority service (PS) policy is a $\mathbf{0}=(0,0)$-standard policy. The PS policy specifies that the server will always choose a customer in (nonempty) queue 2 to
serve at each decision epoch. If queue 2 is empty, the serve will serve customers in queue 1, if there is any. If the server is serving a type 1 customer when a type 2 customer arrives, the type 1 customer is pushed back to queue 1 and the server begins to serve the type 2 customers. The interrupted type 1 customer will resume or repeat its service if the server is available to serve type 1 customers. If the system is empty, the server will be idle.

To prove (i), denote the optimal stationary policy by $\pi^*$. At state $(q_1,q_2)$, we add a virtual customer of type 1 at queue 1. He has the same transfer rate as the ordinary customer of type 1. However, he has no holding cost and transferring cost. For this queueing system $G(q_1,q_2;1,0)$, policy $\pi^*$ can still be used and by comparing each realized trajectory we know that the resulting expected discounted cost $C(G(q_1,q_2;1,0))$ is less than $J_{\alpha}^*(q_1+1,q_2)$. On the other hand, the queueing system $G(q_1,q_2;1,0)$ is in fact a queueing system with state $(q_1,q_2)$ and since policy $\pi^*$ for system $G(q_1,q_2;1,0)$ might not be an optimal policy for queueing system with state $(q_1,q_2)$ we have that $C(G(q_1,q_2;1,0))\geq J^*(q_1,q_2)$. Therefore, $J_{\alpha}^*(q_1+1,q_2)\geq J_{\alpha}^*(q_1,q_2)$ and thus $J_{\alpha}^*(q_1,q_2)$ is increasing in $q_1$. Similarly, $J_{\alpha}^*(q_1,q_2)$ is increasing in $q_2$. Thus, $J_{\alpha}^{*}(\mathbf{q})$ is increasing in $\mathbf{q}$.

Let $\epsilon=\mu_2-\lambda_1-\lambda_2>0$ and $d$ be the PS policy. From \cite{xhz08} it is known that the Markov process induced by $d$ is ergodic, and thus $m_{\mathbf{q},\mathbf{0}}(d)<\infty$, $\forall\ \mathbf{q}\in S$. Next we prove that $c_{\mathbf{q},\mathbf{0}}(d)<\infty$, $\forall\ \mathbf{q}\in S$.

Inspired by \cite{xhz08}, we choose the Lyapunov function $r(\mathbf{q})=Kr_1^{q_1}r_2^{q_2}$ and then apply Lemma~\ref{lem_1} with $H^*=\{\mathbf{0}\}$. Here the constants $K$, $r_1$, $r_2$ are left to be specified later. (\ref{eq_lem_2}) requires that
\begin{eqnarray}
&& c(\mathbf{q})+Kr_1^{q_1}r_2^{q_2}\left[\mu_2\left(\frac{1}{r_2}-1\right)+\lambda_1(r_1-1)+\lambda_2 (r_2-1)+q_1\lambda_T\left(\frac{r_2}{r_1}-1\right)\right]\leq 0, \nonumber \\ && \hspace{9cm}q_1\geq 0,q_2\geq 1,\label{eq_con_1}
\end{eqnarray}
and
\begin{eqnarray}
&&  c(q_1,0)+Kr_1^{q_1}\left[\mu_1\left(\frac{1}{r_1}-1\right)+\lambda_1(r_1-1)+\lambda_2 (r_2-1)+q_1\lambda_T\left(\frac{r_2}{r_1}-1\right)\right]\leq 0, \nonumber
\\ &&\hspace{9cm}  q_1\geq 1.\label{eq_con_2}
\end{eqnarray}

Choose $r_2=r_1$ and $r_1>1$ such that $r_1>\frac{\min(\mu_1,\mu_2)}{\lambda_1+\lambda_2}$. Denote $\delta=\left(\frac{\min(\mu_1,\mu_2)}{r_1}-\lambda_1-\lambda_2\right)\cdot(r_1-1)$. By the choice of $r_1$ it is known that $\delta>0$. Choose $K$ such that $K(r_1-1)>\max(h_1+c\lambda_T, h_2)$. Therefore, we have
\begin{eqnarray*}
&&c(\mathbf{q})+Kr_1^{q_1}r_2^{q_2}\left[\mu_2\left(\frac{1}{r_2}-1\right)+\lambda_1(r_1-1)+\lambda_2 (r_2-1)+q_1\lambda_T(\frac{r_2}{r_1}-1)\right]
\\&\leq& c(\mathbf{q})-Kr_1^{q_1+q_2}\delta
\\&\leq&h_1q_1+h_2q_2+cq_1\lambda_T -K(r_1-1)\delta(q_1+q_2)
\\&\leq& 0,
\end{eqnarray*}
and thus (\ref{eq_con_1}) holds. Similarly, (\ref{eq_con_2}) also holds. Therefore, it follows from Lemma~\ref{lem_1} that $c_{\mathbf{q},\mathbf{0}}(d)<0$ for $\mathbf{q}\neq \mathbf{0}$. Besides, it is easily seen that
$$c_{\mathbf{0},\mathbf{0}}=\frac{\lambda_1}{\lambda_1+\lambda_2}c_{(1,0),\mathbf{0}}(d)+\frac{\lambda_2}{\lambda_1+\lambda_2}c_{(0,1),\mathbf{0}}(d)<\infty.$$
Therefore, the PS policy $d$ is a $\mathbf{0}$-standard policy, and thus (ii) is proved.

It follows from Theorem~\ref{the_2} that an average-cost optimal stationary policy exists and the ACOI is satisfied. Besides, condition (iv) in Theorem~\ref{the_1}(iv) is satisfied as there are only finite $j$ such that $q(j|i,a)\neq 0$ for any $i\in S, a\in A(i)$. Therefore, the ACOE holds for any $ i\in S$.
\end{proof}

{\bf Remark:} From the proof we know that the result still holds if the cost rate function is increasing and polynomial in $\mathbf{q}$.

\section*{Acknowledgements}
The authors gratefully acknowledge that this work was supported by NSFC under grant 71201154, NSFC major program (Grant No. 71090401/71090400) and CPSF under grants 2012M521260 and 2013T60627.


\end{document}